# Matrix Modeling of Energy Hub with Variable Energy Efficiencies

Wujing Huang, *Student Member, IEEE*, Ning Zhang, *Senior Member, IEEE*, Yi Wang, *Member, IEEE*, Tomislav Capuder, *Senior Member, IEEE*, Igor Kuzle, *Senior Member, IEEE*, and Chongqing Kang, *Fellow, IEEE*

*Abstract*—The modeling of multi-energy systems (MES) is the basic task of analyzing energy systems integration. The variable energy efficiencies of the energy conversion and storage components in MES introduce nonlinearity to the model and thus complicate the analysis and optimization of MES. In this paper, we propose a standardized matrix modeling approach to automatically model MES with variable energy efficiencies based on the energy hub (EH) modeling framework. We use piecewise linearization to approximate the variable energy efficiencies; as a result, a component with variable efficiency is equivalent to several parallel components with constant efficiencies. The nonlinear energy conversion and storage relationship in EH can thus be further modeled under a linear modeling framework using matrices. Such matrix modeling approach makes the modeling of an arbitrary EH with nonlinear energy components highly automated by computers. The proposed modeling approach can further facilitate the operation and planning optimization of EH with variable efficiencies. Case studies are presented to show how the nonlinear approximation accuracy and calculation efficiency can be balanced using the proposed model in the optimal operation of EH.

*Index Terms*—multi-energy systems, energy hub, matrix modeling, piecewise linearization, variable efficiency, part-load performance, operation optimization.

## I. INTRODUCTION

The interactions between electricity, gas and heat/cooling systems enable higher energy utilization efficiency and more intermittent renewable energy utilization [1-3]. The modeling of multi-energy systems (MES) has gained increasing importance in recent years as the foundation of the operation and planning of MES, which can unlock the flexibility of shifting across multiple energy vectors and result in reduced costs and lower emissions as compared to separate operation and planning of energy systems. The coupling relationship among different energy sectors in MES and the interaction characteristics between MES and external networks/loads are the key issues that the modeling approach needs to address. The energy hub (EH) concept was introduced to model MES as a unit in which multiple energy carriers can be converted, conditioned and stored. EHs consume energy at their input ports, which are connected to external energy infrastructures, e.g., power grid and gas network. They provide required energy services, such as electricity, heat and cooling, at their output ports [4]. The compactness and effectiveness of EH-based modeling facilitate the calculation of energy flows in MES of different levels, such as trigeneration systems [5], residential consumers [6], commercial buildings [7], industrial parks [8] and national energy systems [9].

While a substantial number of papers have discussed the optimal planning [10] and operation [11] of MES based on the EH model, very few papers have studied the standardized and automatic modeling of MES. Based on the concept of EH, a nonlinear framework for modeling of energy systems with conversion and storage of multiple energy carriers is proposed in [12]. In this modeling framework, the converter and storage coupling matrices are manually formulated according to the given EH layout. Reference [5] proposes an approach to automatically generate coupling matrices for small-scale trigeneration systems using dispatch factors [11]. The formulation involves products of dispatch factors and energy flows, leading to nonlinear rather than linear optimization problems. As a result, the approach is more applicable to small-scale and simple-structure MES. Reference [13] introduces a linear and automatic modeling tool, Hubert, which leverages a concise ASCII description format for EH with an energy "input, input storage, converter, output storage, output" structure. Reference [14] improves Hubert by adding an "export" branch between "converter" block and "output storage" block to facilitate selling PV-generated energy to the external grid. Reference [15] proposes the idea of a linear coupling relationship of an EH by introducing the concept of augmented variables. The standardized procedure for modeling the dispatch-factor-free coupling matrix is then designed to facilitate the computerized calculation for arbitrary EH configurations. Reference [16] proposes an automatic and standardized matrix modeling method for EH with arbitrary configuration. The method uses graph theory to cast the topology and the characteristics of the energy converters into a matrix form and uses linear equations for the dispatch of energy flows without using dispatch factors. However, all the energy efficiencies are regarded as constants in the above research. In reality, energy efficiencies of energy conversion and storage components vary with their operating conditions. The products of variable efficiencies and energy flows bring nonlinearity to the components and systems. Energy conversion efficiencies of converters can reduce significantly when they do not work in rated condition; therefore, the constant efficiency approximation model of MES with high nonlinearity is very likely to be inaccurate and impractical.

There are already some methods taking into consideration the variable efficiencies of energy conversion and storage components. In [17], nonlinear energy converters are directly modeled with highly nonlinear part-load efficiency curves. The resulting model is a nonlinear programming (NP) problem which gives no guarantee that the global optimum can be achieved and few numerically robust solvers for NP problems exist. To avoid NP problems, several linearization approaches are used to model nonlinear converters. The first one is step-wise approximation based on binary variables which guarantees that only one operating point is selected in each time step (e.g., [18] [19]). This formulation carries the disadvantage of only allowing single tabled values to be selected, with no interpolation between them. The second one is piecewise linear (PWL) approximation using special-ordered-set (SOS) variables (e.g., [20]). A special ordered set is a set of consecutive variables in which not more than one (named SOS1) or two (named SOS2) adjacent members may be non-zero in a feasible solution. By introducing SOS2 variables as weighting variables, any operating point can be written as a linear combination of two pre-defined operating points. Although these methods are able to address variable efficiencies, the formulation can hardly be automated by computers. Reference [21] improves Hubert by enabling the modeling of nonlinear converters using PWL approximation. The operating points are defined as the sum of a set of consecutive variables related to different PWL segments, which has been widely used to approximate quadratic fuel cost and transmission loss in power system operation and planning. However, this model can still only deal with the proposed structure of MES where the location of components such as energy storage components are fixed and only one energy conversion block is considered in an EH. In addition, these methods cannot handle the nonlinear energy converters with multiple outputs and complex operation characteristics, e.g., combined heat and power (CHP) in extraction condensing operation mode in which the ratio between electricity and heat production is adjustable within a certain range.

In this paper, we propose a standardized matrix modeling approach for EH with variable energy efficiencies. Using piecewise linearization, the nonlinear energy conversion or charging/discharging process in an energy component is divided into several linear processes with constant efficiencies. Splitters and concentrators are proposed as standardized components to facilitate the split and merge of energy flows imposed by piecewise linearization; as a result, a nonlinear energy component is converted into a "splitter, linear energy component, concentrator" three-layer structure. The nonlinear energy conversion and storage relationship in EH can thus be further modeled under a linear modeling framework using matrices. Such modeling technique facilitates an automatic modeling of an arbitrary EH with nonlinear energy components. Compared to the constant efficiency model and the nonlinear model, the proposed standardized matrix model is proven to significantly improve the accuracy and calculation efficiency, respectively, of the operation optimization of EH. The contributions of this paper, compared to [16], are threefold:

1) Providing a standardized matrix modeling approach for EH with nonlinear energy conversion and storage components, which makes the modeling of an arbitrary EH highly automated by computers. Such technique fills the gap that current nonlinear/linearized modeling approaches can hardly automatically formulate an arbitrary EH with nonlinear energy components.

2) Proposing a generalized linearized model that is capable of accurately modeling any energy conversion or storage component with constant or variable efficiencies under one linear modeling framework. Especially, the generalized linearized model is able to handle the nonlinear energy converters with multiple outputs and complex operation characteristics, e.g., CHP in extraction condensing operation mode.

3) Providing analysis of how the nonlinear approximation accuracy and calculation efficiency can be balanced using the proposed model in the optimal operation of EH.

The rest of this paper is organized as follows: Section II provides the basic backgrounds of the standardized matrix model of EH. Section III proposes the standardized matrix model method to tackle the challenges of variable energy efficiencies of energy components. Section IV conducts a comprehensive case study and applies the proposed modeling method to the operational optimization problem. Section V draws the conclusion.

## II. BACKGROUND: STANDARDIZED MATRIX MODEL OF EH WITH CONSTANT EFFICIENCIES

First, a standardized matrix model of EH with constant energy efficiencies [16] is briefly introduced to facilitate understanding of the proposed matrix model of EH with variable efficiencies and the relationship between these two models.

From the viewpoint of graph theory, each energy conversion or storage component in an EH can be treated as a node, and each energy flow to or from a component can be treated as a branch, as shown in Fig. 1. The inputs and outputs of an EH are treated as special nodes. Each node has several input and output ports through which it exchanges energy with other nodes.

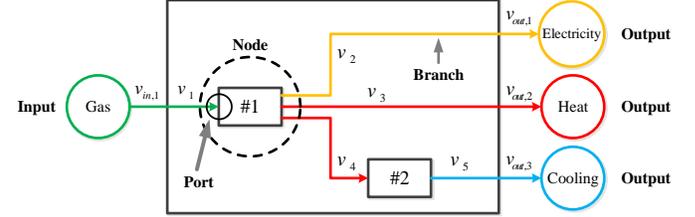

Fig. 1. Definitions of basic elements of an EH in terms of graph theory.

### A. Formulating Port-Branch Incidence Matrices

A port-branch incidence matrix $\mathbf{A}_g$ of node g defines the connections between the ports of node g and all the branches in the EH. If node $g$ has $K_g$ ports and the EH has $B$ branches, the port-branch incidence matrix $\mathbf{A}_g$ has dimensions $K_g \times B$ and is defined as follows:

$$A_g(k,b) = \begin{cases} 1, & \text{branch } b \text{ is connected to input port } k \text{ of node } g \\ -1, & \text{branch } b \text{ is connected to output port } k \text{ of node } g \\ 0, & \text{otherwise} \end{cases} \quad (1)$$

where $k$ is the serial number of input and output ports of a node.

### B. Formulating Converter Characteristic Matrices

The converter characteristic matrix $\mathbf{H}_g$ of node $g$ defines the energy conversion characteristics of the node. If node $g$ has $S_g$ energy conversion processes and $K_g$ ports, the converter characteristic matrix $\mathbf{H}_g$ has dimensions $S_g \times K_g$. If each energy conversion process in node $g$ is related to only one input and one output, then $\mathbf{H}_g$ is defined as follows

$$H_g(s,k) = \begin{cases} \eta_s, & \text{input port } k \text{ is related to} \\ & \text{energy conversion process } s \\ 1, & \text{output port } k \text{ is related to} \\ & \text{energy conversion process } s \\ 0, & \text{otherwise} \end{cases} \quad (2)$$

where $\eta_s$ is the energy conversion efficiency of energy conversion process $s$.

## C. Formulating Energy Conversion Matrices

Given the port-branch incidence matrix $\mathbf{A}_g$ and the converter characteristic matrix $\mathbf{H}_g$, we can calculate the nodal energy conversion matrix $\mathbf{Z}_g$ for node $g$:

$$\mathbf{Z}_g = \mathbf{H}_g \mathbf{A}_g \quad (3)$$

The energy conversion matrix $\mathbf{Z}_g$ for node $g$ defines the relationship between the energy flows in the branches connected to node $g$. This matrix has as many rows as there are energy conversion processes in node $g$ and as many columns as there are branches in the EH.

Suppose that an EH has $N$ nodes. The system energy conversion matrix $Z$ combines the nodal energy conversion matrix of all nodes in the EH:

$$\mathbf{Z} = \begin{bmatrix} \mathbf{Z}_1^{\mathrm{T}} & \mathbf{Z}_2^{\mathrm{T}} & \cdots & \mathbf{Z}_N^{\mathrm{T}} \end{bmatrix}^{\mathrm{T}} \quad (4)$$

The energy flows in all branches form a $B \times 1$ vector $\mathbf{V}$. According to the energy conversion characteristic of each node, we can obtain the energy conversion equation of the EH:

$$\mathbf{ZV} = \mathbf{0} \quad (5)$$

## D. Comprehensive Energy Flow Equations for EH

We define the $m$-dimensional energy input vector $\mathbf{V}_{in}$ and the $n$-dimensional energy output vector $\mathbf{V}_{out}$ if an EH has $m$ inputs and $n$ outputs. The input incidence matrix $\mathbf{X}$ is an $m \times B$ matrix that relates the energy inputs to the branch energy flows. Similarly, the output incidence matrix $\mathbf{Y}$ is an $n \times B$ matrix that relates the energy outputs to the branch energy flows:

$$X(i,j) = \begin{cases} 1, & \text{branch } j \text{ is connected to input port } i \\ 0, & \text{otherwise} \end{cases} \quad (6)$$

$$Y(i,j) = \begin{cases} 1, & \text{branch } j \text{ is connected to output port } i \\ 0, & \text{otherwise} \end{cases} \quad (7)$$

The input incidence and output incidence equations are:

$$\mathbf{V}_{in} = \mathbf{XV}, \quad \mathbf{V}_{out} = \mathbf{YV} \quad (8)$$

The comprehensive energy flow equations of the EH are:

$$\begin{bmatrix} \mathbf{X}^{\mathrm{T}} & \mathbf{Y}^{\mathrm{T}} & \mathbf{Z}^{\mathrm{T}} \end{bmatrix}^{\mathrm{T}} \mathbf{V} = \begin{bmatrix} \mathbf{V}_{in}^{\mathrm{T}} & \mathbf{V}_{out}^{\mathrm{T}} & \mathbf{0}^{\mathrm{T}} \end{bmatrix}^{\mathrm{T}} \quad (9)$$

## III. PROPOSED STANDARDIZED MATRIX MODEL METHODS

When the energy efficiencies of energy conversion and storage components vary with their input or output powers, the constant efficiency model of EH is no longer accurate. Fig. 2 shows the output power curve of a nonlinear converter and its linear approximation which has a significant error at a certain load level. The approximations of variable efficiencies of production units, particularly CHP units, have a high impact on short-term schedules [22] [23]. The products of variable efficiencies and energy flows also bring about nonlinearity to the optimization problem.

In this paper, EH with variable efficiencies is modeled using a piecewise linearization method. Through piecewise linearization, a nonlinear energy conversion or storage component is divided into several equivalent components with constant efficiencies. The modeling of two types of energy converters, i.e., single-input single-output (SISO) and single-input multi-output (SIMO) energy converters, and the modeling of energy storage components are presented. The modeling of energy storage components is similar to that of SISO energy converters where the input port models charging and the output port discharging.

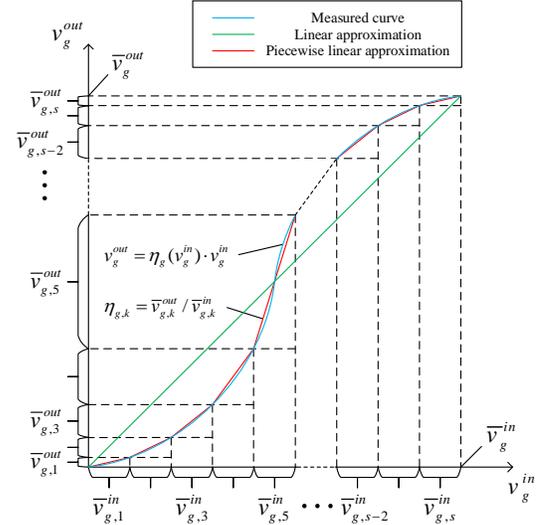

Fig. 2. Nonlinear component variable efficiency linear approximation and piecewise linear approximation.

## A. Piecewise Linearization of Energy Efficiency

Energy conversion efficiency is taken as an example here.

### 1) Type 1: SISO energy converters

Suppose that node $g$ has a single input $v_g^{in}$ and a single output $v_g^{out}$ and that its energy conversion efficiency $\eta_g(v_g^{in})$ is a function of the input energy:

$$v_g^{out} = \eta_g(v_g^{in}) \cdot v_g^{in} \quad (10)$$

As shown in Fig. 2, the range of $v_g^{in}$ (denoted by $\bar{v}_g^{in}$) can be divided into $s$ segments:

$$\begin{aligned} &\left[0, \bar{v}_{g,1}^{in}\right], \left[\bar{v}_{g,1}^{in}, \bar{v}_{g,1}^{in} + \bar{v}_{g,2}^{in}\right], \left[\bar{v}_{g,1}^{in} + \bar{v}_{g,2}^{in}, \bar{v}_{g,1}^{in} + \bar{v}_{g,2}^{in} + \bar{v}_{g,3}^{in}\right], \\ &\cdots, \left[\bar{v}_{g,1}^{in} + \cdots + \bar{v}_{g,s-1}^{in}, \bar{v}_{g,1}^{in} + \cdots + \bar{v}_{g,s}^{in}\right] \end{aligned} \quad (11)$$

where $\bar{v}_{g,k}^{in}$ is the range of the $k$-th segment and $\sum_{k=1}^{s} \bar{v}_{g,k}^{in} = \bar{v}_g^{in}$.

Secondary variables $v_{g,k}^{in} (k = 1, 2, \cdots s)$ are introduced to represent the amounts of input energy that fall in each segment:

$$v_{g,k}^{in} = \begin{cases} \bar{v}_{g,k}^{in}, & v_g^{in} \geq \sum_{i=1}^{k} \bar{v}_{g,i}^{in} \\ v_g^{in} - \sum_{i=1}^{k-1} \bar{v}_{g,i}^{in}, & \sum_{i=1}^{k-1} \bar{v}_{g,i}^{in} < v_g^{in} < \sum_{i=1}^{k} \bar{v}_{g,i}^{in} \\ 0, & v_g^{in} \leq \sum_{i=1}^{k-1} \bar{v}_{g,i}^{in} \end{cases} \quad (12)$$

and $v_{g,k}^{in} (k = 1, 2, \cdots s)$ are subject to

$$v_g^{in} = \sum_{k=1}^{s} v_{g,k}^{in} \quad (13)$$

The range of $v_g^{out}$ can be divided into the same amount of segments as that of $v_g^{in}$:

$$\sum_{i=1}^{k} \bar{v}_{g,i}^{out} = \eta_g \left( \sum_{i=1}^{k} \bar{v}_{g,i}^{in} \right) \cdot \sum_{i=1}^{k} \bar{v}_{g,i}^{in}, \quad k = 1, 2, \cdots, s \quad (14)$$

Similarly, secondary variables $v_{g,k}^{out}$ ($k=1,2,\cdots s$) are introduced and subject to

$$v_g^{out} = \sum_{k=1}^{s} v_{g,k}^{out} \tag{15}$$

In the $k$-th ($k=1,2,\cdots,s$) segment, the variable efficiency $\eta_g(v_g^{in})$ is replaced by a fixed efficiency $\eta_{g,k}$ determined by the secant line of the power output curve in this segment:

$$v_{g,k}^{out} = \eta_{g,k} \cdot v_{g,k}^{in}, \quad k=1,2,\cdots,s \tag{16}$$

$$\eta_{g,k} = \bar{v}_{g,k}^{out} / \bar{v}_{g,k}^{in} \tag{17}$$

Binary variables $u_{g,1}, u_{g,2}, \cdots, u_{g,s-1}$ are introduced to guarantee the continuity of $v_g^{in}$ and $v_g^{out}$:

$$\begin{cases} u_{g,1} \cdot \bar{v}_{g,1}^{in} \leq v_{g,1}^{in} \leq \bar{v}_{g,1}^{in} \\ u_{g,2} \cdot \bar{v}_{g,2}^{in} \leq v_{g,2}^{in} \leq u_{g,1} \cdot \bar{v}_{g,2}^{in} \\ u_{g,3} \cdot \bar{v}_{g,3}^{in} \leq v_{g,3}^{in} \leq u_{g,2} \cdot \bar{v}_{g,3}^{in} \\ \vdots \\ u_{g,s-1} \cdot \bar{v}_{g,s-1}^{in} \leq v_{g,s-1}^{in} \leq u_{g,s-2} \cdot \bar{v}_{g,s-1}^{in} \\ 0 \leq v_{g,s}^{in} \leq u_{g,s-1} \cdot \bar{v}_{g,s}^{in} \end{cases} \tag{18}$$

If $v_g^{in}$ is in the $k$-th ($1<k<s$) segment, then:

$$\begin{cases} u_{g,1}, u_{g,2}, \cdots, u_{g,k-1} = 1 \\ u_{g,k}, u_{g,k+1}, \cdots, u_{g,s-1} = 0 \end{cases} \tag{19}$$

Equation (18) can be written in matrix form to make it easily automated by computers:

$$\bar{\mathbf{V}}_g^{in} \mathbf{U}_{g,a} \leq \mathbf{V}_g^{in} \leq \bar{\mathbf{V}}_g^{in} \mathbf{U}_{g,b} \tag{20}$$

where

$$\begin{cases} \bar{\mathbf{V}}_g^{in} = diag(\bar{v}_{g,1}^{in}, \bar{v}_{g,2}^{in}, \cdots, \bar{v}_{g,s}^{in}) \\ \mathbf{V}_g^{in} = \begin{bmatrix} v_{g,1}^{in} & v_{g,2}^{in} & \cdots & v_{g,s}^{in} \end{bmatrix}^{\mathrm{T}} \\ \mathbf{U}_{g,a} = \begin{bmatrix} \mathbf{U}_g & 0 \end{bmatrix}^{\mathrm{T}}, \ \mathbf{U}_{g,b} = \begin{bmatrix} 1 & \mathbf{U}_g \end{bmatrix}^{\mathrm{T}} \\ \mathbf{U}_g = \begin{bmatrix} u_{g,1} & u_{g,2} & \cdots & u_{g,s-1} \end{bmatrix} \end{cases} \tag{21}$$

*2) Type 2: SIMO energy converters*

If the energy output of each output port is proportional to the energy input, e.g., CHP in backpressure operation mode, the energy conversion processes from the input port to each output port can be piecewise linearized separately in a similar way as SISO energy converters.

A more complicated situation is when the proportion of energy directed to output ports can be flexibly adjusted, e.g., CHP in extraction condensing operation mode. The input-output characteristics become a multivariate quadratic function:

$$F = a \cdot P^2 + b \cdot Q^2 + c \cdot PQ + d \cdot P + e \cdot Q + f \tag{22}$$

where $F$, $P$, $Q$ denote the input gas, output electricity, and output heat, respectively, and $a$, $b$, $c$, $d$, $e$, $f$ are all constants. The characteristics of most SIMO energy converters can be modeled into a multivariate quadratic function, as shown in (22).

We conduct a linear mapping to the input-output characteristics:

$$\begin{cases} F = F_1 + F_2 \\ F_1 = a \cdot \tilde{P}^2 + d \cdot \tilde{P} + f_1 \\ F_2 = \tilde{b} \cdot \tilde{Q}^2 + \tilde{e} \cdot \tilde{Q} + f_2 \end{cases} \tag{23}$$

where:

$$\begin{bmatrix} \tilde{F} \\ \tilde{P} \\ \tilde{Q} \end{bmatrix} = \begin{bmatrix} 1 & 0 & 0 \\ 0 & 1 & c/2a \\ 0 & 0 & 1 \end{bmatrix} \begin{bmatrix} F \\ P \\ Q \end{bmatrix}, \begin{cases} \tilde{b} = b - c^2/4a \\ \tilde{e} = e - cd/2a \\ f = f_1 + f_2 \end{cases} \tag{24}$$

After this mapping, the energy converter can be equivalently seen as two parallel independent ones. Curves $F_1$ and $F_2$ can be separately piecewise linearized in a similar way as SISO energy converters. A broader sense of "output energy" is introduced, i.e., $\tilde{P} = P + (c/2a)Q$, which is a linear combination of the original output energy types. Energy conversion with more output ports can be handled similarly.

### B. Introducing Splitter/Concentrator

Through piecewise linearization, a nonlinear energy conversion or charging/discharging process in a SISO component is divided into several parallel processes with constant efficiencies. From a component level viewpoint, a nonlinear SISO component with input $v_g^{in}$ and output $v_g^{out}$ can be equivalent to a linear multi-input multi-output (MIMO) component with inputs $v_{g,k}^{in}$ ($k=1,2,\cdots,s$) and outputs $v_{g,k}^{out}$ ($k=1,2,\cdots,s$), as shown in Fig. 3.

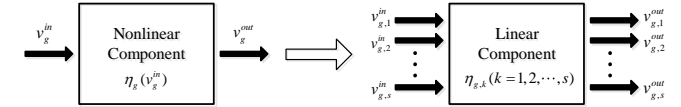

Fig. 3. Piecewise linearization of a nonlinear SISO component.

In this case, the input and output ports of the component should be split into multiple ports. Splitters and concentrators are thus introduced to facilitate the split and merge of branches as shown in Fig. 4. In this paper, the original branches in the EH are named primary branches, and the branches split from primary branches are named secondary branches.

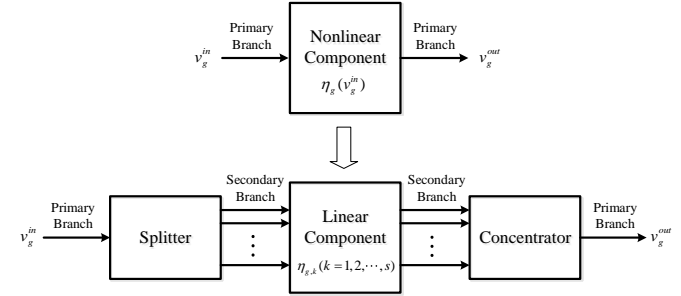

Fig. 4. Splitter and concentrator.

If a nonlinear SIMO converter has one input and $b$ outputs, then one splitter, $b$ concentrators and one linear mapping component should be introduced, as shown in Fig. 5.

The energy flows in all secondary branches form a vector $\mathbf{V}'$ (energy flows in secondary branches connected to the input ports of each node come first in $\mathbf{V}'$), and the primary branch vector $\mathbf{V}$ should be expanded as:

$$\tilde{\mathbf{V}} = \begin{bmatrix} \mathbf{V}^{\mathrm{T}} & \mathbf{V}'^{\mathrm{T}} \end{bmatrix}^{\mathrm{T}} \tag{25}$$

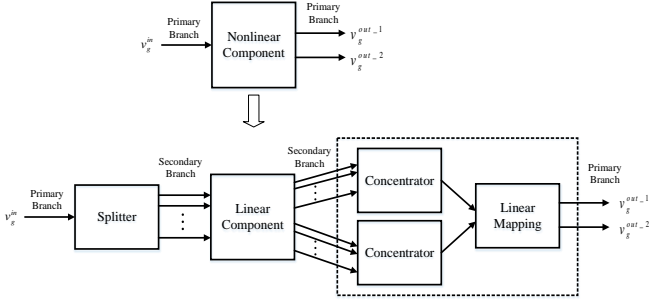

Fig. 5. Piecewise linearization of a nonlinear SIMO energy converter.

## C. Formulating Port-Branch Incidence Matrices

First, consider the simplest case where the EH has only one nonlinear SISO energy converter, which is represented by node $g$, and its input and output ports are connected to only one branch each. Using the method of piecewise linearization, node $g$ is converted into an $s$-input $s$-output converter, and $2s$ secondary branches are introduced. A port-branch incidence matrix $\tilde{\mathbf{A}}_g$ of node $g$ defines the connections between the $2s$ ports and all the branches (including primary and secondary branches). If the EH has $B$ primary branches, then the port-branch incidence matrix $\tilde{\mathbf{A}}_g$ has dimensions $2s \times (B + 2s)$ and is defined as follows:

$$\tilde{\mathbf{A}}_g = \begin{bmatrix} \mathbf{0}_{s \times B} & \mathbf{I}_{s \times s} & \mathbf{0}_{s \times s} \\ \mathbf{0}_{s \times B} & \mathbf{0}_{s \times s} & -\mathbf{I}_{s \times s} \end{bmatrix} \tag{26}$$

where the first $B$ columns correspond to primary branches and the last $2s$ columns correspond to secondary branches.

Then, consider the situation in which the EH has several nonlinear energy converters and several linear energy converters:

1) If node $g$ represents a nonlinear energy converter, $\tilde{\mathbf{A}}_g$ is defined as follows:

$$\tilde{A}_g(k,b) = \begin{cases} 1, & \text{branch } b \text{ is secondary branch and} \\ & \text{is connected to input port } k \text{ of node } g \\ -1, & \text{branch } b \text{ is secondary branch and} \\ & \text{is connected to output port } k \text{ of node } g \\ 0, & \text{otherwise} \end{cases} \tag{27}$$

2) If node $g$ represents a linear energy converter, $\tilde{\mathbf{A}}_g$ is defined as follows:

$$\tilde{\mathbf{A}}_g = \begin{bmatrix} \mathbf{A}_g & \mathbf{0} \end{bmatrix} \tag{28}$$

where $\mathbf{A}_g$ is the port-branch incidence matrix of node $g$ when the EH has no secondary branches and the number of columns of the $\mathbf{0}$ matrix is equal to the length of $\mathbf{V}'$.

## D. Formulating Component Characteristic Matrices

### 1) SISO energy converters

Using the method of piecewise linearization, a nonlinear SISO converter is converted into an $s$-input $s$-output converter. The converter characteristic matrix $\tilde{\mathbf{H}}_g$ has dimensions $s \times 2s$ and is defined as follows:

$$\tilde{\mathbf{H}}_g = \begin{bmatrix} \boldsymbol{\eta}_{\mathbf{g}, s \times s} & \mathbf{I}_{s \times s} \end{bmatrix} \tag{29}$$

where $\boldsymbol{\eta}_{\mathbf{g}, s \times s} = diag(\eta_{g,1}, \eta_{g,2}, \cdots, \eta_{g,s})$ and $\eta_{g,k}$ is the constant energy conversion efficiency in the $k$-th segment.

### 2) SIMO energy converters

Using piecewise linearization, a nonlinear 1-input $b$-output converter is converted into an $s$-input $bs$-output converter. If the energy output of each output port is proportional to the energy input, the converter characteristic matrix $\tilde{\mathbf{H}}_g$ has dimensions $bs \times (b+1) \cdot s$ and is defined as follows:

$$\tilde{\mathbf{H}}_g = \begin{bmatrix} \boldsymbol{\eta}_{\mathbf{1}, s \times s} & \mathbf{I}_{s \times s} & \mathbf{0}_{s \times s} & \cdots & \mathbf{0}_{s \times s} \\ \boldsymbol{\eta}_{\mathbf{2}, s \times s} & \mathbf{0}_{s \times s} & \mathbf{I}_{s \times s} & \ddots & \mathbf{0}_{s \times s} \\ \vdots & \vdots & \ddots & \ddots & \mathbf{0}_{s \times s} \\ \boldsymbol{\eta}_{\mathbf{b}, s \times s} & \mathbf{0}_{s \times s} & \cdots & \mathbf{0}_{s \times s} & \mathbf{I}_{s \times s} \end{bmatrix} \tag{30}$$

If the proportion of energy directed to the output port can be flexibly adjusted, the converter characteristic matrix $\tilde{\mathbf{H}}_g$ has dimensions $s \times (b+1) \cdot s$ and is defined as follows:

$$\tilde{H}_g(i,k) = \begin{cases} 1, & \text{input port } k \text{ is related to energy conversion process } i \\ 1/\eta_i, & \text{output port } k \text{ is related to energy conversion process } i \\ 0, & \text{otherwise} \end{cases} \tag{31}$$

The converter characteristic matrices of linear energy converters remain the same.

## E. Formulating Energy Balance Matrices

Given the port-branch incidence matrix $\tilde{\mathbf{A}}_g$ and the component characteristic matrix $\tilde{\mathbf{H}}_g$, we can calculate the nodal energy balance matrix for node $g$:

$$\tilde{\mathbf{Z}}_g = \tilde{\mathbf{H}}_g \tilde{\mathbf{A}}_g \tag{32}$$

Suppose that EH has $N$ nodes. The system energy balance matrix $\tilde{\mathbf{Z}}$ combines the nodal energy balance matrix of all nodes in the EH:

$$\tilde{\mathbf{Z}} = \begin{bmatrix} \tilde{\mathbf{Z}}_1^{\mathrm{T}} & \tilde{\mathbf{Z}}_2^{\mathrm{T}} & \cdots & \tilde{\mathbf{Z}}_N^{\mathrm{T}} \end{bmatrix}^{\mathrm{T}} \tag{33}$$

The energy balance matrix for each node defines the relationship between the energy flows in the branches connected to the node:

$$\tilde{\mathbf{Z}}\tilde{\mathbf{V}} = \mathbf{0} \tag{34}$$

## F. Formulating Splitter/Concentrator Characteristic Matrices

Splitter/concentrator characteristic matrix $\tilde{\mathbf{W}}_g$ of node $g$ defines the relationship between the primary branches and secondary branches related to node $g$.

First, consider the simplest case where the EH has only one nonlinear SISO energy conversion or storage component, which is represented by node $g$, and its input and output ports are connected to only one branch each. $\overline{\mathbf{A}}_g$ is a matrix whose elements are the absolute value of those in $\mathbf{A}_g$ (the port-branch incidence matrix of node $g$ before piecewise linearization). If the number of segments in piecewise linearization is $s$, then:

$$\overline{\mathbf{A}}_g \mathbf{V} = \mathbf{W}_g \mathbf{V}' \tag{35}$$

$$\mathbf{W}_g = \begin{bmatrix} \mathbf{1}_{1 \times s} & \mathbf{0}_{1 \times s} \\ \mathbf{0}_{1 \times s} & \mathbf{1}_{1 \times s} \end{bmatrix} \tag{36}$$

where $\overline{\mathbf{A}}_g \mathbf{V}$ denotes the energy flows in the primary branches related to node $g$. $\mathbf{1}_{1 \times s}$ ($\mathbf{0}_{1 \times s}$) denotes a $1 \times s$ vector whose elements are all ones (zeros).

We define splitter/concentrator characteristic matrix $\tilde{\mathbf{W}}_g$ of node $g$ as:

$$\tilde{\mathbf{W}}_g = \begin{bmatrix} -\bar{\mathbf{A}}_g & \mathbf{W}_g \end{bmatrix} \quad (37)$$

Then,

$$\tilde{\mathbf{W}}_g \begin{bmatrix} \mathbf{V}^{\mathrm{T}} & \mathbf{V}'^{\mathrm{T}} \end{bmatrix}^{\mathrm{T}} = \tilde{\mathbf{W}}_g \tilde{\mathbf{V}} = \mathbf{0} \quad (38)$$

The splitter/concentrator characteristic matrix for general situations can be formulated in a similar way, which is not explained here due to the length limit.

Suppose that EH has $G$ nonlinear components in the EH. The system splitter/concentrator characteristic matrix $\tilde{\mathbf{W}}$ combines the nodal splitter/concentrator characteristic matrix of all nonlinear components:

$$\tilde{\mathbf{W}} = [\tilde{\mathbf{W}}_1^{\mathrm{T}} \quad \tilde{\mathbf{W}}_2^{\mathrm{T}} \quad \cdots \quad \tilde{\mathbf{W}}_G^{\mathrm{T}}]^{\mathrm{T}} \quad (39)$$

$$\tilde{\mathbf{W}}\tilde{\mathbf{V}} = \mathbf{0} \quad (40)$$

### G. Comprehensive Energy Flow Equations for EH

The input incidence matrix $\tilde{\mathbf{X}}$ is an $m \times \tilde{B}$ ($\tilde{B}$ is the length of $\tilde{\mathbf{V}}$) matrix that relates the energy inputs to the branch energy flows. Similarly, the output incidence matrix $\tilde{\mathbf{Y}}$ is an $n \times \tilde{B}$ matrix that relates the energy outputs to the branch energy flows.

We define:

$$\tilde{\mathbf{X}} = \begin{bmatrix} \mathbf{X} & \mathbf{0} \end{bmatrix}, \quad \tilde{\mathbf{Y}} = \begin{bmatrix} \mathbf{Y} & \mathbf{0} \end{bmatrix} \quad (41)$$

where the number of columns of $\mathbf{0}$ matrix equals to the length of $\mathbf{V}'$.

Then,

$$\mathbf{V}_{in} = \tilde{\mathbf{X}}\tilde{\mathbf{V}}, \quad \mathbf{V}_{out} = \tilde{\mathbf{Y}}\tilde{\mathbf{V}} \quad (42)$$

Based on the above, the comprehensive energy flow equations for EH are:

$$\begin{bmatrix} \tilde{\mathbf{X}}^{\mathrm{T}} & \tilde{\mathbf{Y}}^{\mathrm{T}} & \tilde{\mathbf{Z}}^{\mathrm{T}} & \tilde{\mathbf{W}}^{\mathrm{T}} \end{bmatrix}^{\mathrm{T}} \tilde{\mathbf{V}} = \begin{bmatrix} \mathbf{V}_{in}^{\mathrm{T}} & \mathbf{V}_{out}^{\mathrm{T}} & \mathbf{0}^{\mathrm{T}} & \mathbf{0}^{\mathrm{T}} \end{bmatrix}^{\mathrm{T}} \quad (43)$$

Except for the difference in formulation of the input/output incidence matrix and energy conversion matrix, equations (43) have extra variables for energy flows in secondary branches and equations for splitters/concentrators compared with (9).

An illustrative example based on a simple EH is given in Appendix A to illustrate the standardized matrix modeling method of EH with variable efficiencies.

## IV. CASE STUDY

This section illustrates the application of the proposed standardized matrix modeling method in operational optimization of EH.

### A. System Description

Fig. 6 illustrates an EH which consists of a CHP in backpressure operation mode (heat-to-power ratio is a constant), a compression electric refrigerator group (CERG), an electric heat pump (HP), an auxiliary boiler (AB) and a heat storage component (HS). The energy inputs are electricity and gas, while the outputs consist of cooling, electricity and heat.

The CHP, CERG and HS are considered to exhibit considerable efficiency variations with their input or output powers changing. The electric efficiency of the CHP in partial loads can degrade more than 20% compared to the full-load efficiency [24]. The specific type of CERG in Fig. 6 is an air-cooled centrifugal chiller whose highest coefficient of performance (COP) occurs at a part-load ratio of 0.71–0.84 rather than at full load [25]. The efficiencies of the HP [20] and AB [26] do not differ significantly under part-load conditions; therefore, their efficiencies are assumed to be constant. Table I lists the parameters of the energy components of this system.

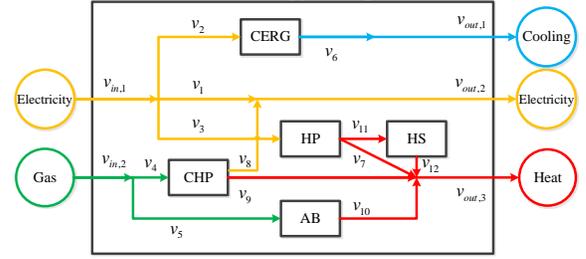

Fig. 6. An EH with three nonlinear energy components.

TABLE I
PARAMETERS OF THE ENERGY COMPONENTS IN FIG. 6

| Converter | Capacity (kW) | Performance [5][20][25][27] |
|---|---|---|
| CERG | 400 | $Out = -0.00003041 \cdot In^3 + 0.01901 \cdot In^2 + 0.2593 \cdot In$ * |
| HP | 400 | $Out = 3 \cdot In$ |
| CHP | El: 300 Therm:420 | El: $Out = 0.0001150 \cdot In^2 + 0.2305 \cdot In$ <br> Therm: $Out = 0.0001611 \cdot In^2 + 0.3228 \cdot In$ ** |
| AB | 900 | $Out = 0.8 \cdot In$ |
| HS | 800 (3.2 MWh) | Charging Efficiency $= -0.00005 * In + 0.93$ <br> Discharging Efficiency $= -0.00005 \cdot Out + 0.93$ |

* Part-load performance of the existing chiller with specific settings of condensing pressure and outdoor temperature; ** Part-load performance of natural gas-powered microturbine with a heat exchanger.

The hourly electricity, heat and cooling demand patterns refer to a hospital site in a Mediterranean area [5] and are indicated in Fig. 7. The gas price is set at 20 euro/MWh and is considered constant during the daily analysis period [5]. The variation in hourly electricity prices for the selected day corresponding to a real case in Italy [5] is also shown in Fig. 7.

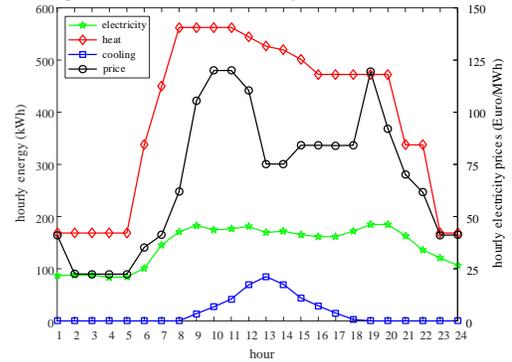

Fig. 7. Hourly load patterns and electricity prices of the selected day.

Using the proposed modeling method described above, the piecewise linearized EH with splitters/concentrators is shown in Fig. 8 (taking 2 segments as an example here).

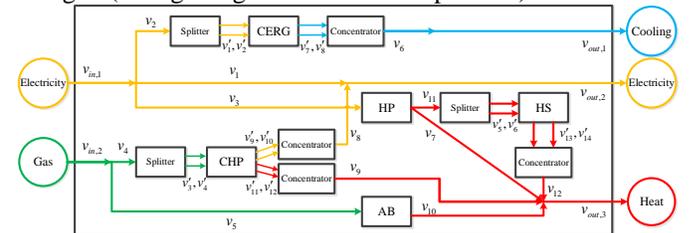

Fig. 8. The piecewise linearized EH with splitters/concentrators.

Given the topology of EH, parameters of energy components and the pre-defined number of segments for piecewise linearization (taking 2 segments as an example here), the computer can automatically formulate the comprehensive energy flow equations of EH as shown in Appendix B.

### B. Multi-Period Optimal Operation Model for EH

The optimal energy flows of a single EH determines how the energy flows should be dispatched among different components to minimize the cost or maximize the profit of the whole system over a given time horizon.

By replacing the comprehensive energy flow equations (9) with (43) and introducing piecewise linearization constraints (e.g., (20)), a revised optimal operation model based on the model proposed in [16] is used.

The decision variables of this problem are the continuous variables $\tilde{\mathbf{V}}$, the continuous variables $\mathbf{V}_{in}$, the continuous variables $\mathbf{E}$ (state of charge of energy storage components) and the binary variables $\mathbf{U}_g$ in each time period, making the problem a mixed integer linear programming (MILP) problem.

### C. Optimization Results and Sensitivity Analysis

Three cases are considered:

1) Case 1: The exact performances of CHP, CERG and HS are considered, making the optimal operation problem a NP problem;

2) Case 2: The efficiencies of the CHP, CERG and HS are approximated as constants (e.g., efficiencies of converters are calculated by dividing maximum output by maximum input), making the optimal operation problem a linear programming problem;

3) Case 3: The proposed standardized matrix modeling model is used. The CHP, CERG and HS are piecewise linearized, making the optimal operation problem an MILP problem. In addition, subcases of Case 3 with different numbers of segments in piecewise linearization are considered. When the number of segments is large enough (e.g., 300), the optimization results are considered accurate and taken as the reference compared with the other optimization results.

Case studies are executed in MATLAB with YALMIP [28] on a PC with an Intel Core i7-7500U 2.70 GHz CPU and 16 GB RAM. Cases 2 and 3 are solved using GUROBI. Case 1 is solved using fmincon with the sequential quadratic programming (SQP) method. The maximum number of iterations is set at 3000. The optimization results of Case 2 are used as initial values for Case 1.

The optimal operating costs, computation times and relative errors of optimal values for these three cases are summarized in Table II and Fig. 9, in which $s$ denotes the number of segments in piecewise linearization. For Case 2, the relative error of the optimal operating cost is over 13%, which indicates that the constant efficiency approximation does not adequately model the EH. With an increase in the number of nonlinear components, nonlinearity of components and complexity of EH, the constant efficiency model of EH results in larger errors due to approximations, making its optimization results impractical. With the proposed standardized matrix modeling method, the relative error of Case 3 is below 1.0% and 0.1% when $s$ is equal to 12 and 36, respectively. With no evident calculation time increase (less than 0.2 seconds) compared to Case 2, the proposed method reduces the approximation error from 13.7% to 0.1% in Case 3. Using the results provided by Case 2 as initial values, Case 1 obtains much more accurate results than Case 2. Nevertheless, the NP problem cannot guarantee a global optimal solution, and its solution is time-consuming. The proposed standardized matrix modeling method makes the optimal operation problem an MILP problem that can be efficiently solved by many existing commercial software (e.g., GUROBI, CPLEX, MOSEK). It should also be noted that the inefficiency of the operation is better approximated when the number of segments in piecewise linearization increases, which leads to a higher operating cost. As shown in Fig. 9, relative error decreases rapidly when $s$ increases from 2 to 20, and computation time increases relatively slowly when $s$ increases from 30 to 70; consequently, the optimal number of $s$ is between 30 and 70 considering the balance between nonlinear approximation accuracy and calculation efficiency.

Table II
OPTIMAL OPERATING COSTS, COMPUTATION TIMES AND RELATIVE ERRORS FOR THE THREE GIVEN CASES

| Case | Optimal Value (Euro) | Computation Time (s) | Relative Error * (%) |
|---|---|---|---|
| Case 1 | 313.628 | 292.73 | 2.35 |
| Case 2 | 264.440 | 0.43 | 13.70 |
| Case 3 ( $s=2$ ) | 280.842 | 0.47 | 8.35 |
| Case 3 ( $s=4$ ) | 289.797 | 0.61 | 5.43 |
| Case 3 ( $s=6$ ) | 297.182 | 0.60 | 3.01 |
| Case 3 ( $s=8$ ) | 301.089 | 0.50 | 1.74 |
| Case 3 ( $s=10$ ) | 303.020 | 0.50 | 1.11 |
| Case 3 ( $s=12$ ) | 304.222 | 0.49 | 0.72 |
| Case 3 ( $s=20$ ) | 305.683 | 0.53 | 0.24 |
| Case 3 ( $s=36$ ) | 306.230 | 0.58 | 0.06 |
| Case 3 ( $s=50$ ) | 306.325 | 0.67 | 0.03 |
| Case 3 ( $s=76$ ) | 306.372 | 0.88 | 0.02 |
| Case 3 ( $s=100$ ) | 306.393 | 1.29 | 0.01 |
| Case 3 ( $s=200$ ) | 306.416 | 2.73 | 0.001 |
| Case 3 ( $s=300$ ) | 306.420 | 3.39 | - |

\* Relative errors for each case/subcase are calculated according to the optimization results of Case 3 ( $s=500$ ).

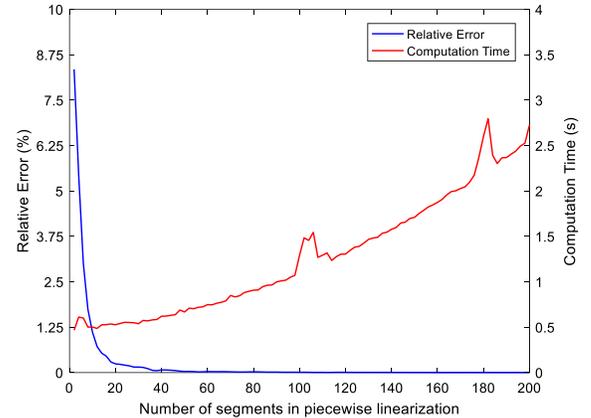

Fig. 9. Computing times and relative errors for different numbers of segments in Case 3.

Fig. 10 shows the operating states of each energy component in the EH over a 24-hour period obtained in Case 2 and Case 3 ( $s=100$ ). The electricity and heat demands are provided by the CHP, HS and AB during the daytime because the price of electricity is high. Load demands and electricity price decrease during the night; therefore, it is inefficient for the CHP to work at a low load level, and it is more cost-effective to purchase electricity from the distribution system to meet the electricity demand and provide power for the HP to produce heat. HS contributes to reducing the heating cost by storing heat when electricity price is low during the night and discharging it in the daytime. The AB is used as a supplemental source for heat supply.

The constant efficiency approximation in Case 2 leads to a wrong and uneconomical operation strategy as shown in Fig. 10: 1) The CHP still works at a relatively low load level during the night without regard to its part-load performance. 2) The electricity input of the CERG decreases because the efficiency is supposed to remain high at any load level. 3) The HS discharges heat with higher power between 10:00 and 14:00 without considering the increasing power loss with increasing discharging power.

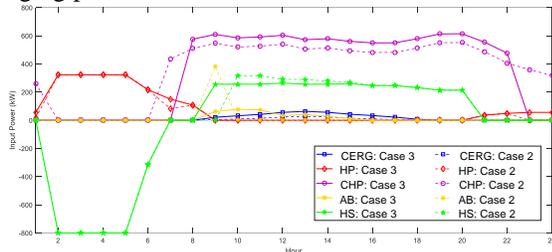

Fig. 10. Operating states of the energy components in the EH of Case 2 and Case 3 ( $s = 100$ ).

## V. CONCLUSION

The variable energy efficiencies of the energy conversion and storage components in MES complicate the modeling of EH in two ways: 1) The nonlinearity caused by variable efficiencies complicates the analysis and optimization of EH; 2) The nonlinear energy components can hardly be automatically modeled by computers. To jointly address these issues, we propose a standardized matrix modeling approach to automatically model EH with variable energy efficiencies. Using piecewise linearization, a nonlinear energy component is converted into a "splitter, linear energy component, concentrator" three-layer structure. The nonlinear energy conversion and storage relationship in EH can thus be further modeled under a linear modeling framework using matrices. Such modeling technique facilitates an automatic modeling of an arbitrary EH with various kinds of nonlinear energy components. Using the proposed modeling approach, the optimal operation of a five-component EH with energy storage shows a great improvement on the accuracy and calculation efficiency compared to the exiting constant efficiency model and nonlinear model, respectively. Future work includes extending the modeling framework to consider the dynamics of energy conversion in different time scales.

## Appendix A: Illustrative Example

An illustrative example based on a simple combined cooling, heating and power (CCHP)-based EH shown in Fig. 1 is given to illustrate the standardized matrix modeling method of EH with variable efficiencies. The EH consists of a combined heat and power (CHP) and a water absorption refrigerator group (WARG).

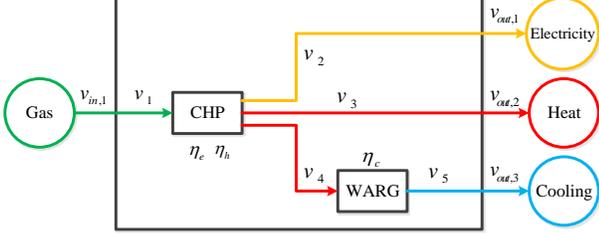

Fig. 1. A CCHP-based EH.

The input and output vectors of this CCHP-based EH are:

$$\mathbf{V}_{in} = \begin{bmatrix} v_{in,1} \end{bmatrix} \quad (1)$$

$$\mathbf{V}_{out} = \begin{bmatrix} v_{out,1} & v_{out,2} & v_{out,3} \end{bmatrix}^{\mathrm{T}} \quad (2)$$

The set of branches (energy flows) can be written as:

$$\mathbf{V} = \begin{bmatrix} v_1 & v_2 & v_3 & v_4 & v_5 \end{bmatrix}^{\mathrm{T}} \quad (3)$$

The CHP is assumed to operate in backpressure mode (the production of electricity and heat are both proportional to the gas input) and produce electricity and heat with fixed efficiencies $\eta_e$ and $\eta_h$, respectively. The WARG is assumed to convert its heat input into cooling with a variable efficiency $\eta_c(v_4)$. Using the piecewise linearization method described in the paper, we can divide the domain of function $\eta_c(v_4)$ into 3 segments (taking 3 segments as an example here) with fixed efficiencies $\eta_{c,1}$, $\eta_{c,2}$, and $\eta_{c,3}$. As a result, the SISO WARG is converted to a 3-input 3-output converter, and one splitter, one concentrator and 6 secondary branches $v_k'(k=1,2,\cdots,6)$ should be introduced. The EH after piecewise linearization is shown in Fig. 2.

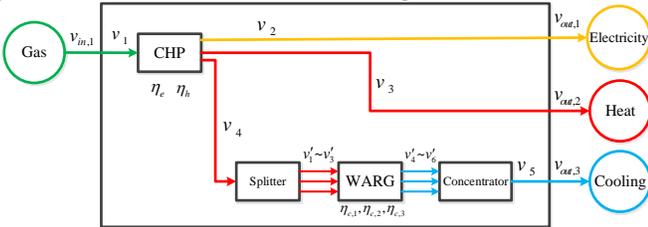

Fig. 2. A CCHP-based EH after piecewise linearization.

The branch vector $\mathbf{V}$ should be expanded:

$$\tilde{\mathbf{V}} = \begin{bmatrix} \mathbf{V} & \mathbf{V}' \end{bmatrix}^{\mathrm{T}} = \begin{bmatrix} v_1 & \cdots & v_5 & v_1' & \cdots & v_6' \end{bmatrix}^{\mathrm{T}} \quad (4)$$

If we assume that the CHP is Node #1 and the WARG is Node #2, the port-branch incidence matrices for the CHP and WARG are:

$$\tilde{\mathbf{A}}_1 = \begin{bmatrix} \mathbf{A}_1 & \mathbf{0} \end{bmatrix} = \begin{bmatrix} 1 & 0 & 0 & 0 & 0 & 0 & 0 & 0 & 0 & 0 & 0 \\ 0 & -1 & 0 & 0 & 0 & 0 & 0 & 0 & 0 & 0 & 0 \\ 0 & 0 & -1 & -1 & 0 & 0 & 0 & 0 & 0 & 0 & 0 \end{bmatrix} \quad (5)$$

$$\tilde{\mathbf{A}}_2 = \begin{bmatrix} \mathbf{0}_{3\times5} & \mathbf{I}_{3\times3} & \mathbf{0}_{3\times3} \\ \mathbf{0}_{3\times5} & \mathbf{0}_{3\times3} & -\mathbf{I}_{3\times3} \end{bmatrix} = \begin{bmatrix} 0 & 0 & 0 & 0 & 0 & 1 & 0 & 0 & 0 & 0 & 0 \\ 0 & 0 & 0 & 0 & 0 & 0 & 1 & 0 & 0 & 0 & 0 \\ 0 & 0 & 0 & 0 & 0 & 0 & 0 & 1 & 0 & 0 & 0 \\ 0 & 0 & 0 & 0 & 0 & 0 & 0 & 0 & -1 & 0 & 0 \\ 0 & 0 & 0 & 0 & 0 & 0 & 0 & 0 & 0 & -1 & 0 \\ 0 & 0 & 0 & 0 & 0 & 0 & 0 & 0 & 0 & 0 & -1 \end{bmatrix} \quad (6)$$

The component characteristic matrices for the CHP and WARG are:

$$\tilde{\mathbf{H}}_1 = \mathbf{H}_1 = \begin{bmatrix} \eta_e & 1 & 0 \\ \eta_h & 0 & 1 \end{bmatrix} \quad (7)$$

$$\tilde{\mathbf{H}}_2 = \begin{bmatrix} \eta_{c,3\times3} & \mathbf{I}_{3\times3} \end{bmatrix} = \begin{bmatrix} \eta_{c,1} & 0 & 0 & 1 & 0 & 0 \\ 0 & \eta_{c,2} & 0 & 0 & 1 & 0 \\ 0 & 0 & \eta_{c,3} & 0 & 0 & 1 \end{bmatrix} \quad (8)$$

Given the port-branch incidence matrices and the component characteristic matrices, we can calculate the nodal energy balance matrix for the CHP and WARG:

$$\tilde{\mathbf{Z}}_1 = \tilde{\mathbf{H}}_1 \tilde{\mathbf{A}}_1 = \begin{bmatrix} \eta_e & -1 & 0 & 0 & 0 & 0 & 0 & 0 & 0 & 0 & 0 \\ \eta_h & 0 & -1 & -1 & 0 & 0 & 0 & 0 & 0 & 0 & 0 \end{bmatrix} \quad (9)$$

$$\tilde{\mathbf{Z}}_2 = \tilde{\mathbf{H}}_2 \tilde{\mathbf{A}}_2 = \begin{bmatrix} 0 & 0 & 0 & 0 & 0 & \eta_{c,1} & 0 & 0 & -1 & 0 & 0 \\ 0 & 0 & 0 & 0 & 0 & 0 & \eta_{c,2} & 0 & 0 & -1 & 0 \\ 0 & 0 & 0 & 0 & 0 & 0 & 0 & \eta_{c,3} & 0 & 0 & -1 \end{bmatrix} \quad (10)$$

The splitter/concentrator characteristic matrix of the WARG is:

$$\tilde{\mathbf{W}}_2 = \begin{bmatrix} -\bar{\mathbf{A}}_2 & \mathbf{W}_2 \end{bmatrix} = \begin{bmatrix} 0 & 0 & 0 & -1 & 0 & 1 & 1 & 1 & 0 & 0 & 0 \\ 0 & 0 & 0 & 0 & -1 & 0 & 0 & 1 & 1 & 1 \end{bmatrix} \quad (11)$$

The input incidence matrix and output incidence matrix of the EH are

$$\tilde{\mathbf{X}} = \begin{bmatrix} \mathbf{X} & \mathbf{0} \end{bmatrix} = \begin{bmatrix} 1 & 0 & 0 & 0 & 0 & 0 & 0 & 0 & 0 & 0 & 0 \end{bmatrix} \quad (12)$$

$$\tilde{\mathbf{Y}} = \begin{bmatrix} \mathbf{Y} & \mathbf{0} \end{bmatrix} = \begin{bmatrix} 0 & 1 & 0 & 0 & 0 & 0 & 0 & 0 & 0 & 0 & 0 \\ 0 & 0 & 1 & 0 & 0 & 0 & 0 & 0 & 0 & 0 & 0 \\ 0 & 0 & 0 & 0 & 1 & 0 & 0 & 0 & 0 & 0 & 0 \end{bmatrix} \quad (13)$$

Based on the above, the comprehensive energy flow equations for the EH are:

$$\begin{bmatrix} 1 & 0 & 0 & 0 & 0 & 0 & 0 & 0 & 0 & 0 & 0 \\ 0 & 1 & 0 & 0 & 0 & 0 & 0 & 0 & 0 & 0 & 0 \\ 0 & 0 & 1 & 0 & 0 & 0 & 0 & 0 & 0 & 0 & 0 \\ 0 & 0 & 0 & 1 & 0 & 0 & 0 & 0 & 0 & 0 & 0 \\ \eta_e & -1 & 0 & 0 & 0 & 0 & 0 & 0 & 0 & 0 & 0 \\ \eta_h & 0 & -1 & -1 & 0 & 0 & 0 & 0 & 0 & 0 & 0 \\ 0 & 0 & 0 & 0 & 0 & \eta_{c,1} & 0 & 0 & -1 & 0 & 0 \\ 0 & 0 & 0 & 0 & 0 & 0 & \eta_{c,2} & 0 & 0 & -1 & 0 \\ 0 & 0 & 0 & 0 & 0 & 0 & 0 & \eta_{c,3} & 0 & 0 & -1 \\ 0 & 0 & 0 & -1 & 0 & 1 & 1 & 1 & 0 & 0 & 0 \\ 0 & 0 & 0 & 0 & -1 & 0 & 0 & 0 & 1 & 1 & 1 \end{bmatrix} \begin{bmatrix} v_1 \\ v_2 \\ v_3 \\ v_4 \\ v_5 \\ v_1' \\ v_2' \\ v_3' \\ v_4' \\ v_5' \\ v_6' \end{bmatrix} = \begin{bmatrix} v_{in,1} \\ v_{out,1} \\ v_{out,2} \\ v_{out,3} \\ 0 \\ 0 \\ 0 \\ 0 \\ 0 \\ 0 \\ 0 \end{bmatrix} \quad (14)$$

Equation (14) provides a matrix model for EH with nonlinear energy conversion components. The proposed matrix modeling approach can be highly automated by computers and facilitates a highly efficient optimization in operation and planning.

APPENDIX B: COMPREHENSIVE ENERGY FLOW EQUATIONS OF EH FOR CASE STUDY

Given the topology of EH, parameters of energy comonents and the pre-defined number of segments for piecewise linearization (taking 2 segments as an example here), the computer can automatically formulate the comprehensive energy flow equations of EH as shown in (15).

$$
\begin{bmatrix}
1 & 1 & 1 & 0 & 0 & 0 & 0 & 0 & 0 & 0 & 0 & 0 & 0 & 0 & 0 & 0 & 0 & 0 & 0 & 0 & 0 & 0 & 0 & 0 & 0 & 0 & 0 \\
0 & 0 & 0 & 1 & 1 & 0 & 0 & 0 & 0 & 0 & 0 & 0 & 0 & 0 & 0 & 0 & 0 & 0 & 0 & 0 & 0 & 0 & 0 & 0 & 0 & 0 & 0 \\
0 & 0 & 0 & 0 & 0 & 1 & 0 & 0 & 0 & 0 & 0 & 0 & 0 & 0 & 0 & 0 & 0 & 0 & 0 & 0 & 0 & 0 & 0 & 0 & 0 & 0 & 0 \\
1 & 0 & 0 & 0 & 0 & 0 & 0 & 1 & 0 & 0 & 0 & 0 & 0 & 0 & 0 & 0 & 0 & 0 & 0 & 0 & 0 & 0 & 0 & 0 & 0 & 0 & 0 \\
0 & 0 & 0 & 0 & 0 & 0 & 1 & 0 & 1 & 1 & 0 & 1 & 0 & 0 & 0 & 0 & 0 & 0 & 0 & 0 & 0 & 0 & 0 & 0 & 0 & 0 & 0 \\
0 & 0 & 0 & 0 & 0 & 0 & 0 & 0 & 0 & 0 & 0 & 0 & 0 & \eta_{CERG,1} & 0 & 0 & 0 & 0 & 0 & 0 & -1 & 0 & 0 & 0 & 0 & 0 & 0 \\
0 & 0 & 0 & 0 & 0 & 0 & 0 & 0 & 0 & 0 & 0 & 0 & 0 & 0 & \eta_{CERG,2} & 0 & 0 & 0 & 0 & -1 & 0 & 0 & 0 & 0 & 0 & 0 & 0 \\
0 & 0 & 3 & 0 & 0 & 0 & -1 & 0 & 0 & 0 & -1 & 0 & 0 & 0 & 0 & \eta^{el}_{CHP,1} & 0 & 0 & 0 & -1 & 0 & 0 & 0 & 0 & 0 & 0 & 0 \\
0 & 0 & 0 & 0 & 0 & 0 & 0 & 0 & 0 & 0 & 0 & 0 & 0 & 0 & 0 & 0 & \eta^{el}_{CHP,2} & 0 & 0 & 0 & -1 & 0 & 0 & 0 & 0 & 0 & 0 \\
0 & 0 & 0 & 0 & 0 & 0 & 0 & 0 & 0 & 0 & 0 & 0 & 0 & 0 & 0 & \eta^{therm}_{CHP,1} & 0 & 0 & 0 & 0 & 0 & -1 & 0 & 0 & 0 & 0 & 0 \\
0 & 0 & 0 & 0 & 0 & 0 & 0 & 0 & 0 & 0 & 0 & 0 & 0 & 0 & 0 & 0 & \eta^{therm}_{CHP,2} & 0 & 0 & 0 & 0 & 0 & -1 & 0 & 0 & 0 & 0 \\
0 & 0 & 0 & 0 & \eta_{HP} & 0 & 0 & 0 & 0 & -1 & 0 & 0 & 0 & 0 & 0 & 0 & 0 & 0 & 0 & 0 & 0 & 0 & 0 & -1 & 0 & 0 & 0 \\
0 & 0 & 0 & 0 & 0 & 0 & 0 & 0 & 0 & 0 & 0 & 0 & -1 & 0 & 0 & 0 & 0 & \eta^{char}_{HS,1} & \eta^{char}_{HS,2} & -1/\eta^{disch}_{HS,1} & -1/\eta^{disch}_{HS,2} & 0 & 0 & 0 & 0 & 0 & 0 \\
0 & 0 & 0 & 1 & 0 & 0 & 0 & 0 & 0 & 0 & 0 & 0 & 0 & -1 & -1 & 0 & 0 & 0 & 0 & 0 & 0 & 0 & 0 & 0 & 0 & 0 & 0 \\
0 & 0 & 0 & 0 & 1 & 0 & 0 & 0 & 0 & 0 & 0 & 0 & 0 & 0 & 0 & -1 & -1 & 0 & 0 & 0 & 0 & 0 & 0 & 0 & 0 & 0 & 0 \\
0 & 0 & 0 & 0 & 0 & 1 & 0 & 0 & 0 & 0 & 0 & 0 & 0 & 0 & 0 & 0 & 0 & -1 & -1 & 0 & 0 & 0 & 0 & 0 & 0 & 0 & 0 \\
0 & 0 & 0 & 0 & 0 & 0 & 1 & 0 & 0 & 0 & 0 & 0 & 0 & 0 & 0 & 0 & 0 & 0 & 0 & -1 & -1 & 0 & 0 & 0 & 0 & 0 & 0 \\
0 & 0 & 0 & 0 & 0 & 0 & 0 & 1 & 0 & 0 & 0 & 0 & 0 & 0 & 0 & 0 & 0 & 0 & 0 & 0 & 0 & -1 & -1 & 0 & 0 & 0 & 0 \\
0 & 0 & 0 & 0 & 0 & 0 & 0 & 0 & 1 & 0 & 0 & 0 & 0 & 0 & 0 & 0 & 0 & 0 & 0 & 0 & 0 & 0 & 0 & -1 & -1 & 0 & 0 \\
0 & 0 & 0 & 0 & 0 & 0 & 0 & 0 & 0 & 1 & 0 & 0 & 0 & 0 & 0 & 0 & 0 & 0 & 0 & 0 & 0 & 0 & 0 & 0 & 0 & -1 & -1
\end{bmatrix}
\begin{bmatrix}
v_1 \\ v_2 \\ v_3 \\ v_4 \\ v_5 \\ v_6 \\ v_7 \\ v_8 \\ v_9 \\ v_{10} \\ v_{11} \\ v_{12} \\ \Delta E_g \\ v'_1 \\ v'_2 \\ v'_3 \\ v'_4 \\ v'_5 \\ v'_6 \\ v'_7 \\ v'_8 \\ v'_9 \\ v'_{10} \\ v'_{11} \\ v'_{12} \\ v'_{13} \\ v'_{14}
\end{bmatrix}
=
\begin{bmatrix}
v_{in,1} \\ v_{in,2} \\ v_{out,1} \\ v_{out,2} \\ v_{out,3} \\ 0 \\ 0 \\ 0 \\ 0 \\ 0 \\ 0 \\ 0 \\ 0 \\ 0 \\ 0 \\ 0 \\ 0 \\ 0 \\ 0 \\ 0
\end{bmatrix}
\tag{15}
$$